\documentclass[11pt,leqno]{article}
\usepackage{graphicx, amsfonts, amsthm, amsxtra, verbatim, multicol}
\usepackage[mathscr]{euscript}
\begin{document}

\begin{center}
{\LARGE\bf On differential modular forms and some analytic relations between
Eisenstein series}\footnote
{

Math. classification: 11F11, 14D07.
}
\\
\vspace{.25in} {\large {\sc Hossein Movasati}} \\
Instituto de Matem\'atica Pura e Aplicada, IMPA\footnote{ Author's address after 1/10/2006.}\\ 
Estrada Dona Castorina, 110 \\
22460-320, Rio de Janeiro, RJ, Brazil\\
 E-mail: {\tt hossein@impa.br} \\
 {\tt http://www.impa.br/$\sim$ hossein/}
\end{center}

\newtheorem{theo}{Theorem}
\newtheorem{exam}{Example}
\newtheorem{coro}{Corollary}
\newtheorem{defi}{Definition}
\newtheorem{prob}{Problem}
\newtheorem{lemm}{Lemma}
\newtheorem{prop}{Proposition}
\newtheorem{rem}{Remark}
\newtheorem{conj}{Conjecture}
\newcommand\diff[1]{\frac{d #1}{dz}} 
\def\End{{\rm End}}              
\def\hol{{\rm Hol}}
\def\sing{{\rm Sing}}            
\def\spec{{\rm Spec}}            
\def\cha{{\rm char}}             
\def\Gal{{\rm Gal}}              
\def\jacob{{\rm jacob}}          
\newcommand\Pn[1]{\mathbb{P}^{#1}}   
\def\Z{\mathbb{Z}}                   
\def\Q{\mathbb{Q}}                   
\def\C{\mathbb{C}}                   
\def\as{\mathbb{U}}                  
\def\ring{{R}}                         
\def\R{\mathbb{R}}                   
\def\N{\mathbb{N}}                   
\def\A{\mathbb{C}}                   
\def\uhp{{\mathbb H}}                
\newcommand\ep[1]{e^{\frac{2\pi i}{#1}}}
\newcommand\HH[2]{H^{#2}(#1)}        
\def\Mat{{\rm Mat}}              
\newcommand{\mat}[4]{
     \begin{pmatrix}
            #1 & #2 \\
            #3 & #4
       \end{pmatrix}
    }                                
\newcommand{\matt}[2]{
     \begin{pmatrix}                 
            #1   \\
            #2
       \end{pmatrix}
    }
\def\ker{{\rm ker}}              
\def\cl{{\rm cl}}                
\def\dR{{\rm dR}}                

\def\hc{{\mathsf H}}                 
\def\Hb{{\cal H}}                    
\def\GL{{\rm GL}}                
\def\pedo{{\cal P}}                  
\def\PP{\tilde{\cal P}}              
\def\cm {{\cal C}}                   
\def\K{{\mathbb K}}                  
\def\k{{\mathsf k}}                  
\def\F{{\cal F}}                     
\def\M{{\cal M}}
\def\RR{{\cal R}}
\newcommand\Hi[1]{\mathbb{P}^{#1}_\infty}
\def\pt{\mathbb{C}[t]}               
\def\W{{\cal W}}                     
\def\Af{{\cal A}}                    
\def\gr{{\rm Gr}}                
\def\Im{{\rm Im}}                
\newcommand\SL[2]{{\rm SL}(#1, #2)}    
\newcommand\PSL[2]{{\rm PSL}(#1, #2)}  
\def\Res{{\rm Res}}              

\def\L{{\cal L}}                     
\def\Aut{{\rm Aut}}              
\def\any{R}                          
\newcommand\ovl[1]{\overline{#1}}    

\def\pm{{\mathsf p}{\mathsf m}}      
\def\T{{\cal T }}                    
\def\tr{{\sf tr}}                 
\newcommand\mf[2]{{M}^{#1}_{#2}}     
\newcommand\bn[2]{\binom{#1}{#2}}    
\def\ja{{\rm j}}                 
\def\Sc{\mathsf{S}}                  
\newcommand\es[1]{g_{#1}}            
\newcommand\V{{\mathsf V}}           
\newcommand\Ss{{\cal O}}             
\def\rank{{\rm rank}}                

\def\Ra{\mathrm{Ra}}
\def\nf{g_2}                         
\begin{abstract}
In the present  article we define the algebra of differential modular forms and we
prove that it is generated by Eisenstein series of weight
$2,4$ and $6$. We define Hecke operators on them, find some analytic relations
between these Eisenstein series and obtain them in a natural way as
coefficients of a family of elliptic curves. 
The fact that a complex manifold over the moduli of polarized
Hodge structures in the case $h^{10}=h^{01}=1$ has an algebraic structure with
an action of an algebraic group plays a basic role in all of the proofs.
\end{abstract}
\section{Introduction}
Around 1970 Griffiths introduced the moduli of polarized Hodge
structures/the period domain $D$ and described
a dream to enlarge $D$ to a moduli space of degenerating polarized
Hodge structures. Since in general $D$ is not a Hermitian symmetric domain,
he asked for the existence of a certain automorphic cohomology theory for
$D$, generalizing the usual notion of automorphic forms on symmetric
Hermitian domains. Since then there have been many efforts in the first
part of Griffiths's dream (see \cite{kaus00, hos006} and the 
references there) but the second part still lives in darkness.

I was looking for some analytic spaces over $D$ for which one may state
Baily-Borel theorem on the unique algebraic structure of quotients
of symmetric Hermitian domains by discrete arithmetic groups. 
I realized that even
in the simplest case of Hodge structures, namely  $h^{01}=h^{10}=1$,  
such spaces are
not well studied. This led me to the definition of a new class of
holomorphic functions 
on the Poincar\'e upper half plane which generalize the classical modular
forms. Since a differential operator acts on them we call them
differential modular forms. These new functions are no longer interpreted
as holomorphic sections of a positive line bundle on some compactified moduli
curve. Nevertheless, they appear in a natural way as coefficients in 
families of elliptic curves, analogous to  Eisenstein series in 
the Weierstrass Uniformization Theorem. 


Recall the Eisenstein series
\begin{equation}
\label{eisenstein}
\es{k}(z)=a_k{\Big (}1+(-1)^k\frac{4k}{B_k}\sum_{n\geq
1}\sigma_{2k-1}(n)e^{2\pi i z n}{\Big )},\ \  k=1,2,3, \ z\in\uhp,
\end{equation}
where $B_k$ is the $k$-th Bernoulli number ($B_1=\frac{1}{6},\
B_2=\frac{1}{30},\ B_3=\frac{1}{42},\ \ldots$), $\sigma_i(n):=
\sum_{d\mid n}d^i$,
\begin{equation}
\label{akhar}
a_1=2\zeta(2)\frac{-1}{2\pi i},\
a_2=2\zeta(4)\frac{60}{(2\pi i)^2},\  a_3=2\zeta(6)\frac{-140}{(2\pi i)^3}
\end{equation}
and $\uhp:=\{x+iy\in \C \mid y>0\}$ is the Poincar\'e upper half plane.
The most well-known differential modular form, which is not a differential
of a  modular form, is the Eisenstein series $\es 1$.
The idea of differentiating modular forms and getting new modular forms
is old and goes back to Ramanujan.
However, the precise definition of differential modular forms has been given
recently in \cite{bu95}. In the present article we give another 
slightly different definition of differential modular forms (see \S \ref{defi}) over a modular subgroup $\Gamma \subset \SL 2\Z$.
It is based on a canonical behavior of holomorphic functions on the
Poincar\'e upper half plane under the action of $\SL 2\Z$. This
approach has the advantage that it can be generalized to any modular
subgroup of $\SL 2\Z$ but the one in \cite{bu95} works only in the case of
full modular group $\SL 2\Z$.
The set of differential modular forms in the present article is a 
bigraded $\C$-algebra
$\mf{}{}=\sum_{n\in\N_0,m\in \N}\mf{n}{m}$, $\mf{0}{m}$ being the set
of classical modular forms of weight $m$, in which the differential
operator $\diff{}$ maps $\mf{n}{m}$ to $\mf{n+1}{m+2}$. We have 
$\es{1}\in \mf{1}{2}, \es{2}\in \mf{0}{4},\es{3}\in\mf{0}{6}$ and we prove:
\begin{theo}
\label{23feb05}
 The functions $g_1,g_2,g_3$ are algebraically
independent and  $\mf{}{}$ is freely generated by
$\es{1},\es{2}$ and $\es{3}$ as a $\C$-algebra.
For  $m\in\N$ and $n\in\N_0$, $\mf nm$ is the set of homogeneous polynomials of degree $m$ in 
the graded ring $\C[g_1,g_2,g_3],\
\deg(g_i)=2i,\ i=1,2,3$ and of degree $2n$ in $g_1$(note that $\deg(g_1)=2$).
\end{theo}
The above theorem implies that $\mf{n}{m}=\{0\}$ for $2n>m$ or $m$ an odd number, 
and  every $f\in \mf{n}{m}$ can be written in a unique way in the form
$\sum_{i=0}^nf_ig_1^i$, where $f_i$ is a modular form of weight $m-2i$.
It generalizes the first theorem in each modular forms book that
the algebra of modular forms is freely generated by the Eisenstein series $g_2$
and $g_3$.  Our proof  gives us also the Ramanujan
relations  between the Eisenstein series $\es{i},\ i=1,2,3$. We define the  action of
Hecke operators on $\mf{n}{m}$ and it turns out that this is similar to the
case of modular forms:
\begin{equation}
\label{quedia}
T_pf(z)=p^{m-n-1} \sum_{d\mid p, 0\leq b\leq d-1}   
d^{-m}f{\Big (}\frac{pz+bd}{d^2}{\Big )},\ p\in \N,\ f\in \mf{n}{m}.
\end{equation}
Hecke operators of this type appear in particular  in the study of
the transfer operator from statistical mechanics which plays an important
role in the theory of dynamical zeta functions (see \cite{HMM}). It turns out that
the differential operator commutes with Hecke operators (see \S \ref{heop}).
Let
$$
g:=(g_1,g_2,g_3): \uhp\rightarrow \C^3
$$
and
$$
T:=\C^3\backslash \{(t_1,t_2,t_3)\in\C^3 \mid 27t_3^2-t_2^3=0\}.
$$
\begin{theo}
\label{realanal}
There are unique analytic functions
$$
B_1, B_2:\ T\rightarrow \R,\  B_3:T\rightarrow \C
$$
such that $B_1$ does not depend on the variable $t_1$ and
\begin{equation}
\label{B1}
B_1\circ g(z)=\Im(z),\ B_1(t_1,t_2k^{-4},t_3k^{-6})=B_1(t)|k|^2
\end{equation}
\begin{equation}
\label{B2}
B_2\circ g=0,\ B_2(
 t_1k^{-2}+k'k^{-1},
t_2k^{-4}, t_3k^{-6})=B_1(t)|k'|^2+B_2(t)|k^{-1}|^2+
\Im(B_3(t)k'\overline{k^{-1}})
\end{equation}
\begin{equation}
\label{B3}
B_3\circ g=1,\ B_3 (
 t_1k^{-2}+k'k^{-1},
t_2k^{-4}, t_3k^{-6})=
B_3(t)k\overline{k^{-1}}+2\sqrt{-1}k
\overline{k'}B_1(t)
\end{equation}
 for all $k\in\C^*$ and $k'\in \C$. Moreover, $|B_3|$ restricted to
 the zero locus of $B_2$ is identically one.
\end{theo}
Differential modular forms are best viewed  as holomorphic functions on
the period domain
\begin{equation}
\label{perdomain}
\pedo:= \left \{
\mat {x_1}{x_2}{x_3}{x_4}\in \GL(2,\C)\mid \Im(x_1\ovl{x_3})>0 \right \},
\end{equation}
so that they are invariant under the action of $\SL 2\Z$ from the left  on $\pedo$ and have 
some compatibility conditions with respect to the action of
\begin{equation}
\label{alggroup}
G_0:=\left \{\mat{k_1}{k_3}{0}{k_2}\mid \ k_3\in\C, k_1,k_2\in \C^*\right \}
\end{equation}
from the right  on $\pedo$ (see Proposition \ref{16aug06}). In this way Theorem \ref{realanal}
is just the translation of the relations of $g_i$'s with there simple analytic functions on $\pedo$ 
(see \ref{chemishod}) into the coefficient space through the period map (see \S \ref{permap}).
The action of Hecke operators on differential modular forms is also best viewed
in this way. We use a four parameter family of elliptic curves in order 
to prove our results on differential modular forms
and in this way we even obtain a result on the periods of the differential forms of the
second type on elliptic curves:
\begin{theo}
\label{ellip}
There is no elliptic curve $E$ and a non-exact differential form of the second type
$\omega$ on $E$, both defined over $\overline{\Q}$, such that
$\int_\delta\omega=0$
for some non-zero topological cycle $\delta\in
H_1(E,\Z)$.
\end{theo}
This theorem   uses Nesterenko's Theorem (see \cite{nes01}) on transcendence
properties of the values of Eisenstein series.
The above theorem for the case in which $\omega$ is of the first kind, 
is well-known. In this case we can even state it for the field $\C$. 
However, it is trivially false when $\omega$ is a differential form of the 
second kind  and we allow transcendental coefficients in $\omega$ or the 
elliptic curve. 

 The present  article stimulates the hope to realize the second
part of Griffiths dream with a new formulation.
The complex manifold $\pedo$  
can  be also introduced over the Griffiths  period domain $D$ 
with an action  of an algebraic group $G_0$
from the right. Since the differential modular forms on $\pedo$ are no longer interpreted 
as sections of positive line
bundles over moduli spaces, the question of the existence of a kind of
Baily-Borel
Theorem for $\pedo$ arises. 
In the case of Hodge structures with $h^{01}=h^{10}=1$ we have $D=\uhp$ and
we show  that $\SL 2\Z\backslash \pedo$ has a canonical
structure of an algebraic quasi-affine variety such that the action of $G_0$
from the right is algebraic.  More precisely, we prove that
$\SL 2\Z\backslash \pedo$ is biholomorphic to
$\C^4\backslash \{t=(t_0,t_1,t_2,t_3)\in \C^4\mid t_0(27t_0t_3^2-t_2^3)=0\}$ and under this
biholomorphism the action of $G_0$ is given by:
\begin{equation}
\label{action}
t\bullet g:=(t_0k_1^{-1}k_2^{-1},
 t_1k_1^{-1}k_2+k_3k_1^{-1},
t_2k_1^{-3}k_2, t_3k_1^{-4}k_2^2), $$ $$
 t=(t_0,t_1,t_2,t_3)\in\C^4,
g=\mat {k_1}{k_3}{0}{k_2}\in G_0.
\end{equation}
The mentioned biholomorphism is given by the period map (see \S \ref{permap}). Using
the methods of present article, one can describe the dynamics of the holomorphic foliation
induced by the Ramanujan's relations. This will be discussed in another paper. 

Let us now explain the structure of this article. 
\S \ref{dmf} is devoted to
the definition of differential modular forms and the action of Hecke operators
on them. 
\S \ref{singular} is devoted to the calculation of the Gauss-Manin
connection of a family of elliptic curves. In this section we
prove that the period map is a biholomorphism and then we take
its  inverse and obtain  the Ramanujan relations.
Finally, \S \ref{proofs} is devoted to 
the proof of theorems  announced in the Introduction.


{\bf Acknowledgment:}
The main ideas of this paper took place in my mind when
I was visiting Prof. Sampei Usui
at Osaka University. Here I would like to thank him for encouraging me to
study Hodge theory and for his help to understand it.
I would like to thank Prof. Karl-Hermann Neeb for  his interest and 
careful reading of the present article. I would like to thank 
the referee of the present article who made useful comments on the first draft of this text and  
introduced me
with the reference \cite{maro05} in which the notion of a differential modular form with the name 
quasi modular form is introduced and Theorem \ref{23feb05} is proved by means of
algebraic methods. The alternative proof presented here by means of the period
map might be useful for further analyzing the differential modular forms and their relations with 
elliptic curves.

\section{$\mf nm$-functions}
\label{dmf}
In this section we use the notations
$A=\mat {a_A}{b_A}{c_A}{d_A}\in\SL 2\R$ and
$$
I=\mat{1}{0}{0}{1},\ T=\mat 1101,\ Q=\mat 0{-1}10,
$$
$$
 x=\mat {x_1}{x_2}{x_3}{x_4},\ 
g=\mat {k_1}{k_2}{0}{k_3},\  x,g\in \GL(2,\C).
$$
When there is no confusion we will simply write  $A=\mat abcd$.
We denote by $\uhp$ the Poincar\'e
upper half plane and
$$\ja (A,z):=c_Az+d_A.$$ For $A\in \SL 2\R$ and $m\in\Z$ we use the
slash operator
$$
f|_mA=(\det A)^{m-1}\ja(A,z)^{-m}f(Az).
$$
For a ring $R$ we denote by $\Mat_p(2,R)$ the set of $2\times 2$-matrices
in $R$ with the determinant $p$.
\subsection{Definitions}
\label{defi}
In this section we define the notion of an $\mf nm$-function. For $n=0$
an $\mf 0m$-function is a classical modular form of weight
$m$ on $\uhp$ (see bellow).
 A holomorphic function $f$ on $\uhp$ is called $\mf nm$
if the following two conditions are satisfied:
\begin{enumerate}
\item
There are holomorphic functions $f_i,\ i=0,1,\ldots,n$ on $\uhp$ such that 
\begin{equation}
\label{4feb05}
f|_m A=\sum_{i=0}^{n}\bn{n}{i} c_A^i\ja(A,z)^{-i}f_i, \ \forall A\in \SL 2\Z.
\end{equation}
\item
$f_i, i=0,1,2,\ldots,n$ have finite growths when $\Im(z)$ tends to $+\infty$, i.e.
$$
\lim_{\Im(z)\to +\infty}f_i(z)=a_{i,\infty} <\infty,\ a_{i,\infty}\in\C.
$$
\end{enumerate}
The above definition can be made using a subgroup $\Gamma\subset \SL 2\Z$. In 
this article we mainly deal with full differential modular forms, 
i.e. the  case $\Gamma=\SL 2\Z$. 
We will also denote  by $\mf nm$ the set of $\mf nm$-functions and
we set
$$
\mf{}{}:=\sum_{m\in\Z, n\in\N_0}\mf nm
$$
For an $f\in\mf{n}{m}$ we have $f|_mI=f_0$ and so $f_0=f$. We have also $f|_mT=f$ and so 
we can write the Fourier expansion of $f$ at
infinity
$$
f=\sum_{n=-N}^{+\infty}a_n q^{n},\ a_n\in \C,\ N=0,1,2,\ldots,\infty,\
\ q=e^{2\pi i z}.
$$
The growth condition on $f$ implies that $N=0$.
Note that for an
$\mf nm$-function $f$ the associated functions $f_i$ are unique.
To see this fix $z$ and consider the right hand side of (\ref{4feb05}) as a polynomial in
$c_A\ja (A,z)^{-1}$ with coefficients $\bn{n}{i}f_{i}$. Since $A$ is an arbitrary element of $\SL 2\Z$
and a one variable polynomial has a finite number of roots, we conclude that $f_i$'s are unique.
\begin{prop}
\label{16august06}
If $f$ is $\mf{n}{m}$-function with the associated functions $f_i$ then $f_i$ is an
$\mf{n-i}{m-2i}$-function with the associated functions $f_{ij}:=f_{i+j},\ j=0,1,\ldots,n-i$, i.e.
\begin{equation}
\label{21feb05}
 f_i|_{m-2i} A=\sum_{j=0}^{n-i}\bn{n-i}{j}
c_A^j\ja(A,z)^{-j}f_{ij}, \ \forall A\in \SL 2\Z,\ f_{ij}=f_{i+j}.
\end{equation}
\end{prop}
\begin{proof}
For $A,B\in \SL 2\Z$ we have
\begin{eqnarray*}
f(ABz) &=&\ja(AB,z)^m\sum_{i=0}^{n}\bn{n}{i} c_{AB}^i\ja(AB,z)^{-i}f_i(z)\\
&=& 
\ja(AB,z)^m\sum_{i=0}^{n}\bn{n}{i} (c_{AB}\ja(B,z))^i \ja(B,z)^{-i}\ja(AB,z)^{-i}f_i(z)\\
&=& 
\ja(AB, z)^m\sum_{i=0}^{n}\sum_{j=0}^i
\bn{n}{i}\bn{i}{j}\ja(AB,z)^jc_{B}^jc_A^{i-j} \ja(B,z)^{-i}\ja(AB,z)^{-i}f_i(z)\\
&=& 
\ja(AB, z)^m\sum_{r=0}^{n}\sum_{s=0}^{n-r}
\bn{n}{r+s}\bn{r+s}{s}c_{B}^sc_A^{r} \ja(B,z)^{-r-s}\ja(AB,z)^{-r}f_{r+s}(z)\\
&=& 
\ja(A, Bz)^m\sum_{r=0}^{n}
\bn{n}{r} c_A^{r} \ja(A,Bz)^{-r}
\left (
\ja(B,z)^{m-2r}
\sum_{s=0}^{n-r}
\bn{n-r}{s}c_{B}^s\ja(B,z)^{-s}
f_{r+s}(z)\right )
\end{eqnarray*}
In the first equality we have used (\ref{4feb05}). 
In the third equality we have used 
$$
\ja(AB,z)c_B+\det(B)c_A=c_{AB}\ja(B,z),\  \forall A,B\in\GL(2,\R).
$$
In the fourth equality we have changed the counting parameters:
$r=i-j,\ s=j,\ 0\leq r+s\leq n$. In the fifth equality we have used
$$
\ja(AB,z)=\ja(A,Bz)\ja(B,z).
$$
From another side
\begin{eqnarray*}
f(ABz)&=& f(A(Bz))\\
&=& \ja(A, Bz)^m\sum_{r=0}^{n}
\bn{n}{r} c_A^{r} \ja(A,Bz)^{-r}f_r(Bz).
\end{eqnarray*}
Since the holomorphic functions associated to $f$ are unique, we conclude that
$$
f_r(Bz)=\ja(B,z)^{m-2r}
\sum_{s=0}^{n-r}
\bn{n-r}{s}c_{B}^s\ja(B,z)^{-s}
f_{r+s}(z),\ \forall B\in\SL2\Z,\ r=0,1,\ldots,n.
$$
\end{proof}
It is useful to define
\begin{equation}
\label{26feb05} f||_mA:=(\det A)^{m-n-1}\sum_{i=0}^{n}\bn{n}{i}
c_{A^{-1}}^i\ja(A,z)^{i-m}f_{i}(Az),\ A\in \GL (2, \R),\ f\in\mf nm.
\end{equation}
The factor $\det A$ is introduced because of Hecke operators
(see \S \ref{heop}).
The
equalities (\ref{4feb05}) is written in the form
\begin{equation}
\label{13feb05} f=f||_{m}A, \forall A\in\SL 2\Z
\end{equation}
(we have substituted $A^{-1}z$ for $z$ and then
$A^{-1}$ for $A$). Since $f||_mA,\ A\in \GL(2,\Z),\ f\in \mf{n}{m}$ may not be in $\mf{}{}$ and it 
is defined using the associated
functions of $f$, it does not make sense to say that $||_m$ is  
an action of $\GL (2,\R)$ on $\mf nm$ from the right. However, we have the following proposition:
\begin{prop}
\label{toulouse}
We have
$$
f||_m A=f||_m(BA),\ \forall A\in \GL(2,\R),\ B\in\SL 2\Z,\ f\in\mf{n}{m}.
$$
\end{prop}
\begin{proof}
The proof is similar to the the proof of Proposition \ref{16august06}.
The term $(\det A)^{n-m+1}f||_mA (z)$ is equal to:
\begin{eqnarray*}
&= & \sum_{i=0}^{n}\bn{n}{i} c_{A^{-1}}^i\ja(A,z)^{i-m}f_i(B^{-1}BAz) \\
& = &
\sum_{i=0}^{n}\sum_{j=0}^{n-i}\bn{n}{i}\bn{n-i}{j}
c_{A^{-1}}^ic_{B^{-1}}^j\ja(A,z)^{i-m}\ja(B^{-1},BAz)^{m-2i-j}f_{i+j}(BAz) \\
&= &
\sum_{r=0}^{n}\sum_{s=0}^{r}
\bn{n}{s}\bn{n-s}{r-s}
c_{A^{-1}}^sc_{B^{-1}}^{r-s}\ja(A,z)^{s-m}\ja(B^{-1},BAz)^{m-r-s}f_{r}(BAz) \\
&= & \sum_{r=0}^{n} \bn{n}{r} \ja(BAz, z)^{r-m}f_r(BAz)
\ja(A,z)^{-r}(\sum_{s=0}^{r}\bn{r}{s}\ja(BA,z)^{s}c_{A^{-1}}^s c_{B^{-1}}^{r-s}) \\
&=& \sum_{r=0}^{n} \bn{n}{r} \ja(BA, z)^{r-m}f_r(BAz)
\ja(A,z)^{-r}(\ja(BA,z)c_{A^{-1}}+c_{B^{-1}})^{r} \\
&=& \sum_{r=0}^{n} \bn{n}{r} \ja(BA,
z)^{r-m}c_{(BA)^{-1}}^rf_r(BAz)= (\det A)^{n-m+1}(f||_mBA)(z)
\end{eqnarray*}
\end{proof}

\subsection{Algebra of $\mf nm$-functions}
\label{algebra}
Recall the Eisenstein series (\ref{eisenstein}) and
\begin{equation}
\label{jdelta}
\Delta(z):=(27g_3^2(z)-g_2^3(z))=-(\frac{2\pi i}{12})^6q\prod_{n=1}^{\infty}
(1-q^n)^{24}=q-24q^2+252q^3+\cdots,
\end{equation}
$$
j(z):=\frac{g_2^3(z)}{-\Delta(z)}=q^{-1}+744+
196884 q+\cdots.
$$
Note that $\zeta(2)=\frac{\pi^2}{6},\zeta(4)=\frac{\pi^4}{90},\zeta(6)=
\frac{\pi^6}{945}$ and so
\begin{equation}
\label{pinfty}
p_\infty:=(a_1,a_2,a_3)=(\frac{2\pi i}{12},12(\frac{2\pi i}{12})^2 ,
8(\frac{2\pi i}{12})^3),
\end{equation}
where $a_i$'s are defined in (\ref{akhar}). 
For $k\geq 2$ one can write
$$
\es {k}(z)=s_k\sum_{0\neq (m,n)\in\Z^2}\frac{1}{(n+mz)^{2k}}\in \mf 0{2k},
$$
where $s_2=\frac{60}{(2\pi i)^2}$ and $s_3=\frac{-140}{(2\pi i)^3}$.
The Eisenstein series $\es{1}$ satisfies
\begin{equation}
\label{g2}
\es{1}\mid_2A-\es{1}=c\ja(A,z)^{-1},\ A\in\SL 2\Z
\end{equation}
and so $\es{1}\in \mf 12$ (see for instance \cite{apo90} p. 69).
The following proposition describes the algebraic structure 
of $\mf {n}{m}$:
\begin{prop}
\label{sarita}
The followings are true:
\begin{enumerate}
\item
For an $f\in \mf 1m$ the function $z(z^{-m}f(\frac{-1}{z})-f(z))$
is in $\mf 0{m-2}$, i.e. it is a modular form of weight $m-2$.
\item
$\mf 12$ is a one dimensional $\C$-vector space generated by $g_1$.
\item
If $n\leq n'$ then $\mf nm\subset \mf {n'}{m}$ and
$$
\mf nm\mf {n'}{m'}\subset \mf {n+n'}{m+m'}, \ \mf nm+\mf {n'}{m}=\mf
{n'}m
$$
\item
For a modular form $f$ of weight $m$ we have $f(g_1)^n\in \mf
n{2n+m}$.
\end{enumerate}
\end{prop}
\begin{proof}
The first item is a direct consequence of Proposition \ref{16august06} and 
the definition of a $\mf{1}{2}$-function applied to $A=Q$:
$$
z^{-m}f(\frac{-1}{z})=f|_mQ=f+z^{-1}f_1(z),\ f_1\in \mf{0}{m-2}.
$$
Since the modular forms of weight $0$ are constant functions, every $f\in \mf{1}{2}$ satisfies:
$f|_2A=f+rc\ja(A,z)^{-1},\ \forall A\in\SL 2\Z$ for some constant $r\in\C$.  This and 
(\ref{g2}) implies that $f-rg_2$ is a modular form of weight $2$. Since there is no non-zero modular
form of weight $2$ we conclude that 
$\mf 12$ is generated by $g_1$.  

If $f\in \mf{n}{m}$ with the associated functions $f_i,\ i=0,1,\ldots,n$ then $f\in \mf{n'}{m}$ with
the associated functions $f_i,\ i=0,1,\ldots,n, \ f_i=0,\ i=n+1,\ldots,n'$. 
If $f\in \mf{n}{m}$ and $g\in \mf{n'}{m'}$ with the associated functions $f_i,\ i=0,1,\ldots,n$
(resp. $g_i,\ i=0,1,\ldots,n'$) then for $A\in\SL 2\Z$
\begin{eqnarray*}
fg|_{m+m'}A&= &f|_mA\cdot g|_{m'}A \\
&=& 
(\sum_{i=0}^{n}\bn{n}{i} c_A^i\ja(A,z)^{-i}f_i)
(\sum_{j=0}^{n'}\bn{n'}{j} c_A^j\ja(A,z)^{-j}g_j)\\
&=& 
\sum_{r=0}^{n+n'} \bn{n+n'}{r} c_A^r\ja(A,z)^{-r}
\left (
\sum_{s=0}^{r}\frac{\bn{n}{s}\bn{n'}{r-s}}{\bn{n+n'}{r}}f_{s}g_{r-s}
\right )
\end{eqnarray*}
which implies that $fg\in \mf{n+n'}{m+m'}$. Now if $m=m'$ then by the discussion at
the beginning of this paragraph we can assume that $n=n'$. Now, $f+g\in\mf{n'}{m}$ with
the associated functions $f_i+g_i,\ i=0,1,\ldots,n'$.

The fourth item  is a consequence of item 3. It
was the main idea behind the definition of $\mf nm$.
\end{proof}
The following proposition shows that $\mf{}{}$ is in fact a
differential algebra.
\begin{prop}
\label{16apr05}
For  $f\in \mf nm$ we have $\diff{f}\in \mf {n+1}{m+2}$ and
\begin{equation}
\label{khordim} \diff{(f||_mA)}= \diff{f}||_{m+2}A,\ \forall  A\in\GL(2,\R).
\end{equation}
\end{prop}
\begin{proof}
For $A\in \GL(2,\Z)$ with $\det(A)=p$ the term $\diff {(f||_mA)}$ is equal to:
\begin{eqnarray*}
 &= &
p^{m-n}\left (\sum_{i=0}^n \bn ni c_{A^{-1}}^i ((m-i)c_{A^{-1}}\ja(A,z)^{i-1-m}
f_i(Az)+\ja(A,z)^{i-m-2}\diff{f_i}(Az))\right ) \\
 &=&
p^{m-n}\left (\sum_{i=1}^{n+1}\bn {n}{i-1}c_{A^{-1}}^i\ja(A,z) ^{i-2-m}(m-i+1)f_{i-1}(Az) \right. \\ 
& + &  
\left. \sum_{i=0}^n \bn ni c_{A^{-1}}^i \ja(A,z)^{i-2-m}\diff{f_i}(Az)\right ) \\
&=&  p^{m-n}\left (\sum_{i=0}^{n+1}\bn {n+1}{i}c_{A^{-1}}^i\ja(A,z) ^{i-2-m}\tilde
f_i(Az)\right ),
\end{eqnarray*}
where
$$\tilde f_i=\frac{i(m-i+1)}{n+1}f_{i-1}+
\frac{n+1-i}{n+1}\diff{f_i},\ i=0,1,\ldots,n+1, \ f_{-1}=f_{n+1}:=0.
$$
For $A\in\SL 2\Z$ we have $f=f||_mA$ and so $\diff{f}=\diff{(f||_mA)}$. This and the above equalities 
imply that $\diff{f}$ is an $\mf{n+1}{m+2}$-function with the associated
$\mf{n-i}{m-2i}$-functions $\tilde f_i,\ i=0,1,2,\ldots,n+1$. The growth
condition on $\tilde f_i$'s follows from
$$
\diff{f}=2\pi i q\frac{df}{dq}
$$
For an arbitrary $A\in\GL (2,\R)$, (\ref{khordim}) follows from the equalities 
at the beginning of the proof.
\end{proof}
The relations between the $g_i,i=1,2,3$ and their derivatives are given by
the Ramanujan's equalities:
\begin{equation}
\label{ramanujan} \diff{g_1}=g_1^2- \frac{1}{12}g_2,\
\diff{g_2}=4g_1g_2- 6g_3,\ \diff{g_3}=6g_1g_3- \frac{1}{3}g_2^2
\end{equation}
(see for instance \cite{la95, nes01}). The proof of Theorem \ref{23feb05} will
contain a new geometric proof of these equalities.
\subsection{Hecke operators}
\label{heop}
For $p\in\N$ let  $\SL 2\Z\backslash  \Mat_p (
2,\Z)=\{[A_1],[A_2],\ldots, [A_s]\}$. We define the $p$-th Hecke
operator in the following way
$$
T_pf:=\sum_{k=1}^s f||_mA_k,\  f\in \mf nm.
$$
Proposition \ref{toulouse} implies that the above definition does not
depend on the choice of $A_k$ in the class $[A_k]$. Form Proposition
\ref{16apr05} one can deduce that the differential operator
$\diff{}$ commutes with the Hecke operator $T_p$.
\begin{prop}
\label{france}
 $T_p$ defines a map from $\mf nm$ to itself.
\end{prop}
This will be proved in \S \ref{heckeproof}.
One can take
$$
\tilde T_p:=\sum_{d\mid p, 0\leq b\leq d-1}\mat {\frac{p}{d}}b0d\in
\Z[\Mat_p(2,\Z)]
$$
and since for matrices $\mat ab0d$ the slash
operator $|_m$ is $p^n$ times $||_m$ we have
$T_pf=p^{-n}f|_m\tilde T_p$
and we get the expression (\ref{quedia}) in the Introduction.
Similar to the case of modular forms (see \cite{apo90} \S 6) one can
check that
$$
T_p\circ T_q=\sum_{d\mid (p,q)}d^{m-n-1}T_{\frac{pq}{d^2}}.
$$
\subsection{The period domain}
The group $\SL 2\Z$ acts from the left on the period domain $\pedo$ defined in
(\ref{perdomain}) and $G_0$ in (\ref{alggroup})
acts  from the right. We consider a holomorphic function on
$$\L:=\SL 2\Z
\backslash \pedo
$$
as a holomorphic function
$$
f:\pedo\rightarrow \C,\hbox{ holomorphic satisfying },
f(Az)=f(z),\ \forall A\in\SL 2\Z, z\in\pedo.
$$
The determinant function  is such a function.
The Poincar\'e upper half plane $\uhp$ is embedded in $\pedo$ in the following
way:
$$
z\rightarrow \tilde z=\mat{z}{-1}{1}{0}.
$$
We denote by $\tilde\uhp$ the image of $\uhp$ under this map. For a
function $f$ on $\uhp$ we denote by $\tilde f$ the corresponding
function on $\tilde \uhp$.
\begin{prop}
\label{16aug06}
There is a unique map  
$$
\phi: \mf{}{}\rightarrow {\cal O}(\pedo),\ f\mapsto \phi(f)=F
$$ of the
algebra of $\mf {}{}$-functions into the algebra of holomorphic functions on
$\pedo$ such that
\begin{enumerate}
\item
For all $f\in \mf {}{}$ the restriction of $F$ to $\tilde\uhp$ is
equal to $\tilde f$.
\item
For all $f\in \mf {}{}$ the holomorphic function $F$ is $\SL 2\Z$
invariant.
\item
We have
\begin{equation}
\label{7feb} F(x\cdot g)=k_2^nk_1^{n-m}\sum _{i=0}^n \bn ni
k_3^ik_2^{-i}F _i(x), \ \forall x\in\pedo,\ g\in G_0,
\end{equation}
where $F_i=\phi(f_i)$.
\end{enumerate}
Conversely, every holomorphic function $F$ on $\pedo$ which
is left $\SL 2\Z$-invariant and satisfies (\ref{7feb}) for some holomorphic functions $F_i$ on $\pedo$ 
such that the restriction of $F_i$'s to $\tilde \uhp$ have  finite
growths at infinity is of the form $F=\phi(f)$ for some $f\in \mf
nm$.
\end{prop}
\begin{proof}
We have
\begin{equation}
\label{neeb}
\mat{x_1}{x_2}{x_3}{x_4}=\mat{\frac{x_1}{x_3}}{-1}{1}{0}
\mat{x_3}{x_4}{0}{\frac{\det(x)}{x_3}}.
\end{equation}
Therefore, we expect $F$ to be defined by
\begin{equation}
\label{khoshhal}
F(x)=F\left ( \mat{\frac{x_1}{x_3}}{-1}{1}{0}
\mat{x_3}{x_4}{0}{\frac{\det(x)}{x_3}}\right ):= x_3^{-m}\det(x)^n \sum
_{i=0}^n \bn ni x_4^ix_3^i\det(x)^{-i}f_i(\frac{x_1}{x_3}).
\end{equation}
Let us prove that the function $f\mapsto F=\phi(f)$ satisfies the items 1,2 and 3.
For $x=\tilde z,\ z\in\uhp$ we have $x_4=0$ and so $F(x)=f_0(z)=f(z)$. 
This proves the first item.

 By the definition of $F$ one
can rewrite (\ref{4feb05}) in the form
\begin{equation}
\label{17aug06}
f(A\frac{x_1}{x_3})=(cx_1+dx_3)^{m-n}F\mat{x_1}{-d}{x_3}{c},
\end{equation}
where $A=\mat{a}{b}{c}{d}\in\SL 2\Z$.

Now, we prove item 3. Let
$$
g'=\mat {k_1'}{k_3'}{0}{k_2'}:=\mat{x_3}{x_4}{0}{\frac{\det(x)}{x_3}}.
$$
\begin{eqnarray*}
\hbox{ RHS of (\ref{7feb}) } &= &
 k_2^nk_1^{n-m}\sum _{i=0}^n
\bn ni k_3^ik_2^{-i}F_i(x) \\
&= &
(k_2k_2')^n(k_1k_1')^{n-m}\sum _{i=0}^n\sum _{j=0}^{n-i}
\bn ni \bn {n-i}j k_3^ik_2^{-i}
k_2'^{-i}k_1'^{i}
k_3'^jk_2'^{-j}f_{i+j}(\frac{x_1}{x_3}) \\
& =&
(k_2k_2')^n(k_1k_1')^{n-m}\sum _{r=0}^n\sum _{s=0}^{r}
\bn ns \bn {n-s}{r-s} k_3^sk_2^{-s}
k_2'^{-s}k_1'^{s}
k_3'^{r-s}k_2'^{-r+s}f_{r}(\frac{x_1}{x_3}) \\
&=&
(k_2k_2')^n(k_1k_1')^{n-m}\sum _{r=0}^n
\bn nr  
\left (
\sum _{s=0}^{r}
\bn{r}{s}(k_2k_3')^{r-s}(k_3k_1')^{s}
\right )
(k_2k_2')^{-r}f_r(\frac{x_1}{x_3}) \\
&=& 
(k_2k_2')^n(k_1k_1')^{n-m}\sum _{r=0}^n
\bn nr  (k_2k_3'+k_3k_1')^r(k_2k_2')^{-r}f_r(\frac{x_1}{x_3}) \\
& = &
F(\mat{z}{-1}{1}{0}g'g)=F(xg) \\
\end{eqnarray*}
In the second equality we have used the definition of $F_i=\phi(f_i)$ 
as in  (\ref{khoshhal}) and
the fact that the associated functions of $f_i$ are $f_{i+j},\ j=0,1,\ldots,n-i$, where
$f_i$'s are associated functions of $f$ (see Proposition \ref{16august06}).
In the third equality we have changed the counting parameters:  $s=i, r=i+j$.
The fifth equality is just the expansion of $(k_2k_3'+k_3k_1')^r$. The sixth equality
is by the definition of $F=\phi(f)$. 

Let us now prove the second item.
We have to prove that $F(Ax)=F(x),\ \forall A\in \SL 2\Z$: The term $F(Ax)$ is equal to 
\begin{eqnarray*}
 &= &
(cx_1+dx_3)^{-m}\det(x)^n
\sum _{i=0}^n \bn ni (cx_2+dx_4)^i(cx_1+dx_3)^i\det(x)^{-i}
f_i(A\frac{x_1}{x_3}) \\
& =&
(cx_1+dx_3)^{-m}\det(x)^n \\ & . &
\sum _{i=0}^n \bn ni (cx_2+dx_4)^i(cx_1+dx_3)^i\det(x)^{-i}
(cx_1+dx_3)^{m-2i-(n-i)}F_i\mat{x_1}{-d}{x_3}{c} \\
& =&
F\left (\mat{x_1}{-d}{x_3}{c}
\mat{1}{\frac{cx_2+dx_4}{cx_1+dx_3}}{0}{\frac{\det(x)}{cx_1+dx_3}}\right )=F(x).
\end{eqnarray*}
In the second equality we have used the fact that $f_i\in\mf{n-i}{m-2i}$ and the corresponding
equality (\ref{17aug06}). In the third equality we have used  the third item of Proposition \ref{16aug06}.

We have finished the proof of the fact that $F=\phi(f)$ has the desired properties.
Now, let $F$ satisfy 2,3 and its restriction to $\tilde \uhp$ has
a finite growth at infinity.  Put $f=F\mid_{\tilde \uhp}$ and $f_i:=F_i\mid_{\tilde \uhp}$. 
We are going to prove that $f$ satisfies (\ref{4feb05}) with the associated functions $f_i$'s  
and so $f\in\mf{n}{m}$. 
First, we note
that
$$
\mat abcd
\mat{z}{-1}{1}{0}=\mat{Az}{-1}{1}{0}\mat{\ja(A,z)}{-c}{0}{\ja(A,z)^{-1}\det(A)},\
A\in\GL (2,\R).
$$
Now
\begin{eqnarray*}
f(Az) &= & F\mat{Az}{-1}{1}{0} =F(\mat abcd
\mat{z}{-1}{1}{0}\mat{\det(A)\ja(A,z)^{-1}}{c}{0}{\ja(A,z)}) \\
 &=& \ja(A,z)^{n}\ja(A,z)^{m-n}\sum _{i=0}^n \bn ni
 c^i\ja(A,z)^{-i}f_i(x)=\sum _{i=0}^n \bn ni
 c^i\ja(A,z)^{m-i}f_i(x).
\end{eqnarray*}
In the second equality we have used the facts that $F$ is $\SL 2\Z$ invariant and
it satisfies the property (\ref{7feb}).  
We have finished  the proof of our proposition.
\end{proof}
We denote by $\check {\mf nm}$ the set of holomorphic functions on  $\pedo$
which  are left $\SL 2\Z$ invariant and 
satisfy (\ref{7feb}) for some holomorphic functions $F_i$ on $\pedo$
such that $F_i$'s restricted to $\tilde \uhp$ have finite growths at infinity.
For the determinant function $\det:\pedo\rightarrow \C$ we have:
\begin{equation}
\label{refreekhub}
\det \in \check{\mf 10},\  (\det)^i\check{\mf{n}{m}}\subset 
\check{\mf{n+i}{m}},\ i\in\N_0.
\end{equation}
In a similar way as in Proposition \ref{sarita}, one can prove that:
\begin{equation}
\label{refreekhub2}
\check{\mf nm}\check{\mf {n'}{m'}}\subset \check{\mf {n+n'}{m+m'}}, \ \check{\mf nm}+(\det)^{n'-n}\check{ \mf {n'}{m}}=
\check{
\mf{n'}
{m}
}
\ m,m'\in\Z,\ n,n'\in\N_0, n\leq n'.
\end{equation}
We have an isomorphism $\check{\mf{n}{m}}\rightarrow \mf{n}{m},\ F\mapsto F|_{\tilde \uhp}$
whose inverse is given by $\phi$ in Proposition \ref{16aug06}.  
For a classical modular form $f:\uhp\rightarrow \C$ of weight $m$ the
associated $F=\phi(f)$ is
$$
F(x)=x_3^m f(\frac{x_1}{x_3})\in \check{\mf 0m}.
$$
\subsection{Proof of Proposition \ref{france}}
\label{heckeproof}
We define
\begin{eqnarray*}
\check T_p:\check {\mf
nm}\rightarrow \check{\mf nm}, \ 
\check T_pF(x)=p^{m-2n-1}
\sum_{k=1}^s F(A_ix).
\end{eqnarray*}
This function has trivially its image in $\check{\mf nm}$. We
calculate the corresponding function in $\mf nm$: The term $T_pf(z)$ is equal to:
\begin{eqnarray*}
 &= & p^{m-2n-1} \sum_{k=1}^s F(A_k\mat{z}{-1}{1}{0})
\\
& = &
p^{m-2n-1}\sum_{k=1}^s
 F\left (\mat{A_kz}{-1}{1}{0}\mat{\ja(A_k,z)}{-c}{0}{p.\ja(A_k^{-1},A_kz)}
\right ) \\
&= &
p^{m-2n-1}\sum_{k=1}^s (p.\ja(A_k^{-1},A_kz))^n(\ja(A_k,z))^{n-m}\sum
_{i=0}^n \bn ni
(-cp^{-1})^i\ja(A_k^{-1},A_kz)^{-i}f_i(A_kz) \\
&=&\sum_{k=1}^s  f||_m A_k,
\end{eqnarray*}
where $A_k=\mat{a}{b}{c}{d}$. This proves Proposition \ref{france}.
\subsection{Some non-holomorphic functions on the period domain}
\label{chemishod}
We now define some functions that will be used in the proof of Theorem 
\ref{realanal}. Their relation with  the functions of Theorem \ref{realanal} will be
explained in the next sections, where we have introduced the period map.

On the complex manifold $\pedo$ we have the following left $\SL 2\Z$ invariant
analytic functions:
$$
B_1:=\Im (x_1\overline{x_3}),\ B_2:=\Im (x_2\overline{x_4}),\
B_3:=x_1\overline{x_4}-x_2\overline{x_3}.
$$
They define analytic functions on $\L$ which we denote them by the same letter.
They satisfy
\begin{equation}
\label{shamlu1}
B_1\mid_{\tilde \uhp}(z)=\Im(z),\ B_1(xg)=B_1(x)|k_1|^2
\end{equation}
\begin{equation}
\label{shamlu2}
B_2\mid_{\tilde \uhp}(z)=0,\ B_2(xg)=B_1(x)|k_3|^2+B_2(x)|k_2|^2+
\Im(B_3(x)k_3\overline{k_2})
\end{equation}
\begin{equation}
\label{shamlu3}
B_3\mid_{\tilde \uhp}(z)=1,\ B_3(xg)=B_3(x)k_1\overline{k_2}+2\sqrt{-1}k_1
\overline{k_3}B_1(x).
\end{equation}
By the equality (\ref{neeb}) one can easily see that
every point in $\pedo$ can be mapped to a point of
$\tilde \uhp$ by an action of a unique element of $G_0$.
This implies that  the $\SL 2\Z$ invariant functions $B_i,\ i=1,2,3$, 
with the above properties are unique.
\section{Families of elliptic curves and the Gauss-Manin connection}
\label{singular}
In this section we consider the following family of elliptic curves
\begin{equation}
\label{4mar05}
 E: \ y^2-4t_0(x-t_1)^3+t_2(x-t_1)+t_3=0
\end{equation}
and its specialization  $E_t$ over a regular point $t\in T:=\C^4\backslash \{\Delta=0\}$,
where $\Delta=t_0(27t_0t_3^2-t_2^3)$  is the discriminant of $E$. 
We have a proper smooth morphism $E\rightarrow T$ defined over $\C$. 
The family $E$ with $t_0=1$ and $t_1=0$ is the classical Weierstrass family of
elliptic curves and the material of the sections \ref{gmcon}, \ref{permap} and \ref{mainproof}
for such a family is well-known (see for instance \cite{ka73} Appendix 1,\cite{gri} and \cite{sas}).
The discussions related to the full family (\ref{4mar05}) are slight modifications of the classical ones. 

\subsection{Gauss-Manin connection} 
\label{gmcon}
The algebraic definition of the Gauss-Manin connection is made by N. M. Katz, T. Oda 1968 and 
P. Deligne 1971. Its computational aspects are discussed in \cite{hos005}. Let us
introduce the basic notations for the  Gauss-Manin connection of the family $E$.

The Gauss-Manin connection on the cohomology bundle ${\cal H}_\dR^1(E/T)$ is a 
$\C$-linear map:
$$
\nabla:{\cal H}_\dR^1(E/T) \rightarrow \Omega^1_{T}\otimes_{{\cal O}_T} {\cal H}_\dR^1(E/T),
$$
where $\Omega^1_T$ is the sheaf of differential $1$-forms on $T$. 
It satisfies $\nabla(fe)=df\otimes e+ f\nabla(e)$, $f$ (resp. $e$) being a section of ${\cal O}_T$ 
(resp. ${\cal H}_\dR^1(E/T)$). The set 
$$
H:=H_\dR^1(E/T)
$$ 
of global sections of ${\cal H}_\dR^1(E/T)$ is a $\C[t,\frac{1}{\Delta}]$-module
generated freely by the classical differential forms $\frac{dx}{y},\frac{xdx}{y}$.
For $\omega=(\frac{dx}{y}, \frac{xdx}{y})^{\tr}$, the Gauss-Manin connection 
can be written in the following way:
\begin{equation}
\label{violette}
\nabla\omega= A\otimes\omega,\
A=\frac{1}{\Delta}(\sum_{i=0}^3 A_i dt_i),\
A_i\in\Mat (2,\C[t]).
\end{equation}
A simple calculation shows that:
{\tiny
\begin{equation}
\label{rosa} A_0 =\mat {3/2t_0t_1t_2t_3-9t_0t_3^2+1/4t_2^3}
{-3/2t_0t_2t_3}
{3/2t_0t_1^2t_2t_3+9t_0t_1t_3^2-1/2t_1t_2^3+1/8t_2^2t_3}
{-3/2t_0t_1t_2t_3-18t_0t_3^2+3/4t_2^3}
\end{equation}
$$
A_1=\mat 0 0 {27t_0^2t_3^2-t_0t_2^3} 0
$$
$$
A_2 =\mat {-9/2t_0^2t_1t_3+1/4t_0t_2^2} {9/2t_0^2t_3}
{-9/2t_0^2t_1^2t_3+1/2t_0t_1t_2^2-3/8t_0t_2t_3}
{9/2t_0^2t_1t_3-1/4t_0t_2^2}
$$
$$
A_3 =\mat {3t_0^2t_1t_2-9/2t_0^2t_3} {-3t_0^2t_2}
{3t_0^2t_1^2t_2-9t_0^2t_1t_3+1/4t_0t_2^2}
{-3t_0^2t_1t_2+9/2t_0^2t_3}.
$$
}
(See \cite{hos005} for the procedures which calculate all matrices above). 
Let $U$ be an small open set in $U$ and $\{\delta_t\}_{t\in U},
\delta_t\in H_1(E_t,\Z)$ be a continuous family of
topological one dimensional cycles. The main property of the Gauss-Manin
connection is:
\begin{equation}
\label{khodayakomak} d(\int_{\delta_t}\eta)=\sum \alpha_i
\int_{\delta_t}\beta_i,\ \nabla\eta=\sum_{i}\alpha_i\otimes
\beta_i,\ \alpha_i\in H^0(T,\Omega^1_{T}),\ \eta, \beta_i\in H.
\end{equation}

\subsection{Period map}
\label{permap}
The period
map associated to the basis $\omega:=(\frac{dx}{y}, \frac{xdx}{y})^{\tr}$ is given by:
$$
\pm: T\rightarrow \SL 2\Z\backslash \pedo,\ t\mapsto
\left [\frac{1}{\sqrt{2\pi i}}\mat
{\int_{\delta_1}\omega_1}
{\int_{\delta_1}\omega_2}
{\int_{\delta_2}\omega_1}
{\int_{\delta_2}\omega_2} \right ].
$$
It is well-defined and holomorphic. Here $\sqrt{i}=e^{\frac{2\pi
i}{4}}$ and $(\delta_1,\delta_2)$ is a basis of the $\Z$-module 
$H_1(E_t,\Z)$
such that the intersection matrix in this basis is
$\mat{0}{1}{-1}{0}$. It follows from (\ref{khodayakomak}) that
$\pm$ satisfies the differential equation:
\begin{equation}
\label{4mar}
 d(\pm)(t)=\pm(t) \cdot A ^\tr,\ t\in T,
\end{equation}
where $d$ is the differential map.
\subsection{The Action of an algebraic group}
\label{mainproof}
We consider the family of elliptic curves (\ref{4mar05}).
It can be checked easily that (\ref{action}) is an action of $G_0$ on $\A^4$
(this can be also  verified from the proof of the proposition bellow).
It is also easy to verify that $\A^4/G_0$ is isomorphic to $\Pn 1$ through
the map
\begin{equation}
s: \A^4/G_0\rightarrow \Pn 1,\  t\rightarrow [t_2^3:27t_0t_3^2-t_2^3]
\end{equation}
and so
\begin{equation}
\label{jinv}
j(t):=\frac{t_2^3}{27t_0t_3^2-t_2^3}
\end{equation}
is $G_0$-invariant and gives an isomorphy between $T/G_0$ and $\A$.
\begin{prop}
\label{cano} The period $\pm$ associated to the basis $\omega$  is a
biholomorphism and
\begin{equation}
\label{gavril}
\pm(t\bullet g)=\pm(t)\cdot g,\ t\in\A^4,\ g\in G_0.
\end{equation}
\end{prop}
\begin{proof}
We first prove (\ref{gavril}). Let
$$
\alpha: \A^2\rightarrow \A^2,\ (x,y)\mapsto
(k_2^{-1}k_1x-k_3k_2^{-1}, k_2^{-1}k_1^{2}y).
$$
Then
$$
k_2^{2}k_1^{-4}\alpha^{-1}(f)=y^2-4t_0k_2^{2}k_1^{-4} (
k_2^{-1}k_1x-k_3k_2^{-1}-t_1)^3+ t_2k_2^{2}k_1^{-4}(
k_2^{-1}k_1x-k_3k_2^{-1}-t_1)+t_3k_2^{2}k_1^{-4}
$$
$$
y^2-4t_0k_1^{-1}k_2^{-1}(x-(t_1k_2k_1^{-1}+k_3k_1^{-1}))^3+
t_2k_1^{-3}k_2(x-(t_1k_2k_1^{-1}+k_3k_1^{-1}))+t_3k_1^{-4}k_2^{2}
$$
This implies that $\alpha$ induces an isomorphism of elliptic curves
$$
\alpha: E_{t\bullet g}\rightarrow E_t.
$$
Now
$$
\alpha^{-1} 
\omega= \mat{k_1^{-1}}{0}{-k_3k_2^{-1}k_1^{-1}}{k_2^{-1}}\omega=\mat {k_1}{0}{k_3}{k_2}^{-1} \omega
$$
and so
$$
\pm(t)= \pm(t\bullet g).g^{-1}
$$
which proves (\ref{gavril}).

Let $B$ be the $4\times 4$ matrix 
whose  $i$-th row, $i=1,2,\ldots,4$, constitutes of the first and
second rows of $A_{i-1}$. A simple calculation shows that
$$
\det(B)=\frac{3}{4}t_0\Delta^3
$$
and so the period map $\pm$ is regular at each point $t\in
T$. Therefore,  it is locally a biholomorphism.

The period map $\pm$ induces a local biholomorphic map
$\bar \pm: T/G_0\rightarrow \SL 2\Z\backslash \uhp\cong \C$ and so we have the
local biholomorphism $\bar \pm \circ j^{-1} : \A\rightarrow \A$.
One can compactify $\SL 2\Z\backslash \uhp$ by adding the cusp $\SL
2\Z/\Q=\{c\}$ (see \cite{la95}) and the map $ \bar \pm \circ j^{-1}$ is
continuous at $v$ sending $v$ to $c$, where $v$ is the point induced by
$27t_0t_3^2-t_2^3=0$ in $\A^4/G_0$.
 Using Picard's Great Theorem we conclude that $j^{-1}\circ \bar \pm$ is a
biholomorphism and so $\pm$ is a biholomorphism.
\end{proof}
\subsection{The inverse of the period map}
We denote by
$$
F=(F_0,F_1,F_2,F_3):\pedo \rightarrow T
$$
the composition of the quotient map $\pedo\rightarrow \SL 2\Z\backslash \pedo$  and  the inverse of the period map.
\begin{prop}
\label{F0123}
The following is true:
\begin{enumerate}
\item
$F_0(x)=\det(x)^{-1}$.
\item
For $i=2,3$
$$
F_i=\det(x)^{1-i}  \check{\es{ i}}\in \check{\mf 0{2i}}
$$
where $\es {i}$ is the Eisenstein series (\ref{eisenstein}).
\item
$F_1=\check{g_1} \in \check {\mf {1}{2}}$.
\end{enumerate}
\end{prop}
\begin{proof}
Taking $F$ of (\ref{gavril}) we have
$$
F_0(xg)=F_0(x)k_1^{-1}k_2^{-1}, 
$$
\begin{equation}
\label{2apr05}
F_1(xg)=F_1(x)k_1^{-1}k_2+k_3k_1^{-1},\
\end{equation}
$$
F_2(xg)=F_2(x)k_1^{-3}k_2,\  F_3(xg)=F_3(x)k_1^{-4}k_2^2,\ \forall 
x\in\L,\
g\in G_0.
$$
By the Legendre's  theorem $\det(x)$ is equal to one on $V:=\pm 
(1\times 0\times \A\times \A)$ and so the same is true for
$F_0\det(x)$. But the last function is invariant under the action of
$G_0$ and so it is the constant function $1$.  This proves the first
item.
Let $G_i=F_i\det(x)^{i-1},\ i=1,2,3$. The equalities 
(\ref{2apr05}) imply that $G_i, i=2,3$ do not depend on $x_2,x_4$.
Now the map $(t_2,t_3)\rightarrow \pi\circ \pm(1,0,t_2,t_3)$, where $\pi$ is 
the projection on the $x_1,x_3$ coordinates, is the classical period map
(see for instance the appendix of 
\cite{ka73}) and this implies that 
$G_i= \check {\es{i}},\ i=2,3$. 
Note that in our definition of the 
period map the factor $\frac{1}{\sqrt{2\pi i}}$ appears.
In particular $F_i, i=2,3$  have finite growths at infinity. 
The fact that $F_1$ has a finite growth at infinity follows form
the Ramanujan relations (\ref{ramanujan}) and the equality 
$\frac{d}{dz}=2\pi i q \frac{d}{dq}$.
Since $G_1\in \check {\mf 12}$, $\check {\mf 12}$ is
a  one dimensional
space, both $g_1,G_1$ satisfy (\ref{2apr05})
and ${\mf 0{2}}=\{0\}$, we have $G_1=g_1$. 
\end{proof}
\subsection{Ramanujan  relations}
We proved in Proposition \ref{cano} that the period map $\pm$ associated
to $\omega$ is a biholomorphism. According to (\ref{4mar}), the
inverse $F$ of $\pm$ satisfies the differential equation
$$
x.A(F(x))^\tr=I.
$$
We consider $\pm$ as a map sending the vector $(t_0,t_1,t_2,t_3)$ to
$(x_1,x_2,x_3,x_4)$. Its derivative at $t$ is a $4\times 4$ matrix whose
$i$-th column constitutes of the first and second row of
$\frac{1}{\Delta}xA_i^\tr$. We  use (\ref{rosa})
to derive the equality {\tiny
$$
(dF)_x=(d\pm)_t^{-1}=
$$
$$
\det(x)^{-1}\left ( \begin{array}{cccc} -F_0x_4 &F_0x_3
&F_0x_2 &-F_0x_1
\\ \frac{1}{12F_0}(12F_0F_1^2x_3-12F_0F_1x_4-F_2x_3)
&-F_1x_3+x_4 &\frac{1}{12F_0}(-12F_0F_1^2x_1+12F_0F_1x_2+F_2x_1)
&F_1x_1-x_2
\\4F_1F_2x_3-3F_2x_4-6F_3x_3
&-F_2x_3 &-4F_1F_2x_1+3F_2x_2+6F_3x_1 &F_2x_1
\\\frac{1}{3F_0}(18F_0F_1F_3x_3-12F_0F_3x_4-F_2^2x_3)
&-2F_3x_3 &\frac{1}{3F_0}(-18F_0F_1F_3x_1+12F_0F_3x_2+F_2^2x_1)
&2F_3x_1
\end{array} \right ).
$$
} 
For $g_i:=F_i\mid_{\tilde \uhp}$ the first column of the above equality 
gives us the Ramanujan
relations (\ref{ramanujan}).
\subsection{The family $y^2-4t_0x^3+t_1x^2+t_2x+t_3$}
The family (\ref{4mar05}) can be rewritten in the
form
$$
y^2-4t_0x^3+12t_0t_1x^2+(-12t_0t_1^2+t_2)x+(4t_0t_1^3-t_2t_1+t_3)=0.
$$
The mapping $$ \alpha: \A^4\rightarrow \A^4, t\mapsto (t_0,
12t_0t_1, -12t_0t_1^2+t_2,  4t_0t_1^3-t_2t_1+t_3) $$ is an
isomorphism and so we can restate Proposition \ref{cano} for the family
$$
E_t: y^2-4t_0x^3+t_1x^2+t_2x+t_3=0.
$$
The inverse of the period map in this case is given by
$G=(G_0,G_1,G_2,G_3)$ with
$$
G_0=F_0,\ 
G_1=12F_0F_1,\ 
G_2=-12F_0F_1^2+F_2,\ 
G_3=4F_0F_1^3-F_2F_1+F_3.
$$
In this case the singular fibers are parameterized by the zeros of
$$
\Delta:=
t_0
(432t_0^2t_3^2+72t_0t_1t_2t_3-16t_0t_2^3+4t_1^3t_3-t_1^2t_2^2).
$$
The Ramanujan relations take the simpler form:
\begin{equation}
\label{ramanu}
\left \{ \begin{array}{l}
\dot {t_1}= -t_2 \\
\dot{t_2}= -6t_3  \\
\dot {t_3}= t_1t_3-\frac{1}{4}t_2^2
\end{array} \right.,
\end{equation}
where
$$
(t_1,t_2,t_3):=(12g_1, -12g_1^2+g_2,4g_1^3-g_2g_1+g_3).
$$
\section{Proofs}
\label{proofs}
Now we are in a position to prove the theorems announced in the Introduction.
\subsection{Proof of Theorem \ref{23feb05}}
It is enough to prove that 
$$
\check{\mf {}{}}:=\sum_{m, i\in\Z,n\in \N_0} 
F_0^{i}\check{\mf{n}{m}}
$$
 as a $\C(F_0)$-algebra is freely generated by $F_i,\ i=1,2,3$ and  every 
 $F\in\check{\mf nm}$ 
 can be written as a homogeneous polynomials of degree $m$ in 
the graded ring $\C(F_0)[F_1,F_2,F_3],\
\deg(F_i)=2i,\ i=1,2,3$ and of degree $2n$ in $F_1$.
  Since $F_0=(\det)^{-1}\mid_{\tilde \uhp}=1$ and $F_i\mid_{\tilde \uhp}=\tilde 
g_i,\ i=1,2,3$, this will imply Theorem \ref{23feb05}.

Since the period map is a biholomorphism (Proposition \ref{cano}) 
and the pull-back of $F_i,\ i=0,1,2,3$ by the period map is $t_i$ and $t_1,t_2,t_3$ are algebraically
independent over $\C(t_0)$,  we conclude that 
$F_i,\ i=1,2,3$ are algebraically independent over $\C(F_0)$. 



Again we use that fact that the period map is a biholomorphism and conclude that 
for $\tilde F\in\check{\mf{n}{m}}$ and 
its associated functions $\tilde F_i\in \check{\mf {n-i}{m-2i}}$, 
there exist holomorphic functions 
$
p_i: T\to \C,\  i=0,1,\ldots,n,\ p_0:=p
$ 
such that  $\tilde F_i=p_i(F_0,F_1,F_2,F_3), \ i=0,1,2,\ldots,n$.
The property (\ref{7feb}) of $\tilde F$ and (\ref{gavril}) imply that:
\begin{equation}
\label{jabr}
p(t\bullet g)=k_{2}^nk_1^{n-m}\sum_{i=0}^n\bn nik_3^ik_2^{-i}p_i(t),\ 
\forall g\in G_0,\ t\in T.
\end{equation}
Take $g=\mat{1}{t_1}{0}{1}$ and $t=(t_0,0,t_1,t_3)$. Then
\begin{equation}
\label{visa}
p(t_0,t_1,t_2,t_3)=\sum_{i=0}^n\bn nit_1^ip_i(t_0,0,t_2,t_3).
\end{equation}
This implies that $p$ is a polynomial of degree at most $2n$ in 
the variable $t_1$ ($\deg(t_1)=2$). 
In (\ref{jabr}) we take $g=\mat{k_1}{0}{0}{t_0k_1^{-1}}$ and obtain
\begin{equation}
\label{baha}
p(1,t_1t_0k_1^{-2}, t_2t_0k_1^{-4},t_3t_0^2k_1^{-6})=t_0^nk_1^{-m}p(t).
\end{equation}
We substitute (\ref{visa}) in (\ref{baha}) and
consider the equalities obtained by the coefficients of 
$t_1^i$. We get 
\begin{equation}
\label{13apr}
p_i(1,0,t_2t_0k_1^{-4}, t_3t_0^2k_1^{-6})=
t_0^{n-i}k_1^{-m+2i}p_i(t_0,0,t_2,t_3),\ i=0,1,2,\ldots,n. 
\end{equation}
We take $t_0=1$ and conclude that $p_i(1,0,F_2,F_3)\in\check{\mf{0}{m-2i}}$. 
Since every modular form of weight $m-2i$ can be written as a homogeneous polynomial
of degree $m-2i$ in $\C[g_2,g_3],\ \deg(g_2)=4,\deg(g_3)=6$, $p_i(1,0,F_2,F_3)$ can be written
as a homogeneous polynomial of degree $m-2i$ in $\C[F_2,F_3],\ \deg(F_2)=4,\ \deg(F_3)=6$.
In (\ref{13apr}) we put $k_1=1$ and conclude
that 
$p_i(t_0,0,t_2,t_3)=t_0^{i-n} p_i(1,0,t_2t_0, t_3t_0^2)$ is a homogeneous polynomial
of degree $m-2i$ in  $\C(t_0)[t_2,t_3],\ \deg(t_2)=4,\deg(t_3)=6$.

\subsection{Proof of Theorem \ref{realanal}}
In \S \ref{chemishod} we described some analytic functions $B_i,\ i=1,2,3$,
on $\L$ which  have the compatibility properties (\ref{shamlu1}),(\ref{shamlu2}) and 
(\ref{shamlu3}) with the action of $G_0$ on
$\L$. We use Proposition \ref{cano} and transfer them to the world of
coefficients $T$. We obtain analytic functions $B_1,B_2:T\to \R$ and $B_3:T\to \C$ which satisfy:
\begin{equation}
\label{shaml1}
 B_1\circ g(z)=\Im(z),\ B_1(t\bullet g)=B_1(t)|k_1|^2
\end{equation}
\begin{equation}
\label{shaml2}
B_2\circ g(z)=0,\ B_2(t\bullet g)=B_1(t)|k_3|^2+B_2(t)|k_2|^2+
\Im(B_3(t)k_3\overline{k_2})
\end{equation}
\begin{equation}
\label{shaml3}
B_3\bullet g(z)=1,\ B_3(t\bullet g)=B_3(t)k_1\overline{k_2}+2\sqrt{-1}k_1
\overline{k_3}B_1(t).
\end{equation}
for $t\in T,\ z\in\uhp$ and $g\in G_0$, where $g:\uhp\to T,\ g(z)=(1,g_1(z),g_2(z),g_3(z))$ 
(for the sake of simplicity we have used the same letters to name these
functions).
In (\ref{shaml1}) we put $g=\mat{1}{-t_1}{0}{1}$ and obtain $B_1(t_0,0,t_2,t_3)=B_1(t_0,t_1,t_2,t_3)$
which means that $B_1$ does not depend on $t_1$.  Now in (\ref{shaml1}) we put 
$t_0=1$ and $g=\mat{k}{0}{0}{k^{-1}}$ and obtain (\ref{B1}).
In (\ref{shaml2}) and (\ref{shaml3}) we put $t_0=1$ and $g=\mat{k}{k'}{0}{k^{-1}}$ and obtain the equalities 
(\ref{B2}) and (\ref{B3}).
The uniqueness of $B_1,B_2$ and $B_3$ in Theorem \ref{realanal} follows form uniqueness of the 
same functions in the period domain.

The proof of the last part of the theorem is as follows: 
If $B_2(x):=\Im(x_2\bar x_4)=0$ for some $x\in\pedo$ with $\det(x)=1$ then 
$x=\mat
{x_1}{x_4r}{x_3}{x_4}$ for some $r\in \R$ and $x_4(x_1-rx_3)=1$. Then
\begin{equation}
\label{10apr05}
B_3(x)=\overline{x_4}(x_1-rx_3)=\frac{\overline{x_4}}{x_4}.
\end{equation}
which implies that $|B_3(x)|=1$. In the coefficient space $\det(x)=1$ corresponds to $t_0=1$ and
so we get the last statement of Theorem \ref{realanal}.

One can say something more about $B_1$:
The function
$
B_1\cdot |\Delta|^\frac{1}{6}
$
is $G_0$ invariant and so there is an analytic function
$b_1:\A\rightarrow \R$
such that
$$
B_1(t)=\frac{b_1(j(t))}{|\Delta(t)|^\frac{1}{6}}.
$$
Taking this equality to the period domain and restricting it to $\tilde \uhp$, we get
$$
\Im(z)=\frac{b_2(j(z))}{|\Delta(z)|^\frac{1}{6}}
$$
where the above $j$ and $\Delta$ are the ones on \S \ref{algebra}.

\subsection{Proof of Theorem \ref{ellip}}
Let $k$ be an algebraically closed field of charachteristic $0$ ,
for instance take $k=\overline{\Q}$.  By a variety over $k$ we mean the set 
of its $k$-rational points. We redefine:
$$
G_0:=\left \{\mat{k_1}{k_3}{0}{k_2}\mid \ k_3\in k, k_1,k_2\in k^*\right \},\ 
T:=k^4\backslash\{t\in k^4\mid t_0(27t_0t_3^2-t_2^3)=0\}.
$$
\begin{prop}
\label{classspace}
The quasi affine variety
$T$ is the moduli of $(F, [\omega_1],[\omega_2])$'s, where $F$ is an elliptic
curve defined over $k$, $\omega_1$ is a differential form of the first kind 
on $F$ and $([\omega_1], [\omega_2])$ is a basis of $H^1_\dR(F)$.
\end{prop}
\begin{proof}
For simplicity we do not write more $[.]$ for differential forms.  
The $j$ invariant (\ref{jinv}) classifies the ellipric curves over $k$ (see 
\cite{har77} Theorem 4.1). Therefore, for a given elliptic curve $F/k$ we can
find parameter a $t\in k^4$ such that $F\cong E_t$ over $k$. 
Under this isomorphism we write
$$
\matt{\omega_1}{\omega_2}=g^\tr\matt{\frac{dx}{y}}{\frac{xdx}{y}},\ 
\hbox{ in }\ H^1_\dR(E_t)
$$
for some $g\in G_0$, where $\omega_1,\omega_2$ are as in the proposition. 
Now, the triple $(F,\omega_1,\omega_2)$ is isomorphic to 
$(E_{t\bullet g}, \frac{dx}{y}, \frac{xdx}{y})$.
Since $j: \A^4/G_0\rightarrow \A$ is an isomorphism, every triple
$(F,\omega_1,\omega_2)$ is represented exactly by one parameter $t\in T$. 
\end{proof}
Let us take $k=\bar\Q$.
By Proposition \ref{classspace} the hypothesis of Theorem \ref{ellip} gives 
us a parameter $t\in T$ such
that $\int_{\delta}\frac{xdx}{y}=0$, for some $\delta\in H_1(
E_t,\Z)$. We can assume that $\delta$ is not a multiple of another cycle
in $H_1(E_t,\Z)$ and so we can find another cycle $\delta'$ in $H_1(E_t,\Z)$ such that
$\langle \delta,\delta'\rangle=1$.  
The corresponding period matrix $x$ of $E_t$ in 
$(\delta',\delta)$ has zero $x_4$-coordinate
and so the numbers
$$
t_0=\det(x)^{-1},\
t_i=F_i(x)=\det(x)^{1-i}x_3^{-2i}\es{i}(\frac{x_1}{x_3}),\ i=2,3,\
t_1=F_1(x)=\det(x)x_3^{-2}\es{1}(\frac{x_1}{x_3})
$$
all are in $\overline{\Q}$. Here we have used Proposition \ref{F0123} and $x_4=0$.
Now, for $z=\frac{x_1}{x_3}\in\uhp$ we have
$$
\frac{\es{3}}{\es{1}^3}(z), \frac{\es{2}}{\es{1}^2}(z),\
\frac{\es{3}^2}{\es{2}^3}(z)\in \overline{\Q}.
$$
This is in contradiction with

{\bf Theorem} (Nesterenko 1996, \cite{nes01})
{\it For any $z\in\uhp$, the set
$$
e^{2\pi i z},\ \frac{g_1(z)}{a_1}, \frac{g_2(z)}{a_2},\frac{g_3(z)}{a_3}
$$
contains at least three algebraically independent numbers over $\Q$. }

A direct corollary of Theorem \ref{ellip} is that the multi-valued
function
$$
I(t)=\frac{\int_{\delta_t}\frac{xdx}{y}}{
\int_{\delta_t}\frac{dx}{y} }
$$
defined in $T$ never takes algebraic values for algebraic $t$.
\subsection{Other topics}
The literature of modular forms and its applications in number theory
is huge. The first question which naturally arises at this point is as
follows: Which part of the theory of modular forms can be generalized to
the context of differential modular forms and which arithmetic properties
can one expect to find? Since I am not expert in this area, I just mention
some subjects  which could fit well into this section.


One may ask for the Eichler-Manin-Shimura theory  of periods for cusp
forms (see \cite{HMM} and its references) in the context of differential modular forms.
Note that the notion ``period'' in this theory, as far as I know,
has nothing to do with the notion of a period in this article. The notion
of period appears there because classical modular forms can be interpreted
as sections of a tensor product of the cotangent bundle of a moduli curve and
hence a differential multi form, which can be integrated over some path in
the moduli curve (see \cite{sho80}). 
The differential modular forms are no longer interpreted
as sections of line bundles and this makes the situation more difficult.
Lewis type equations attached to differential modular forms will  be also 
of interest (see \cite{HMM}).

Another theory which could be developed for differential modular forms is
Atkin-Lehner theory of old and new modular forms (see the references in \cite{HMM}).
This seems to me to be a quite accessible theory. The $L$-functions attached
to differential modular forms through their Fourier expansion and the extension of
the Rankin-Cohen bracket to differential modular forms may be also  of
interest.


\def\cprime{$'$} \def\cprime{$'$}

\bibliographystyle{plain}

\end{document}